%% file: UZV_arXiv.tex
\documentclass[journal,12pt,onecolumn,draftclsnofoot,]{IEEEtran}

\ifCLASSINFOpdf
\else
\fi

\usepackage{cite}
\usepackage{amsmath}
\usepackage{amssymb}
\newtheorem{theorem}{Theorem}

\usepackage{footnote}
\makesavenoteenv{tabular}
\usepackage{algorithm}
\usepackage{algorithmic}

\usepackage{mleftright}
\usepackage{epstopdf}
\usepackage{caption}
\usepackage{subcaption}

\usepackage{tikz}
\usetikzlibrary{shapes,arrows,fit,positioning,shadows,calc}
\usetikzlibrary{plotmarks}
\usetikzlibrary{decorations.pathreplacing}
\usetikzlibrary{patterns}
\usetikzlibrary{automata}
\usepackage{pgfplots}
\pgfplotsset{compat=newest}
\usepackage{adjustbox}
\usepackage{caption}
\begin{document}
\title{Randomized Rank-Revealing UZV Decomposition for Low-Rank Approximation of Matrices}

\author{\IEEEauthorblockN{Maboud F. Kaloorazi, and 
Rodrigo C. de Lamare} \\
\vspace{5mm}
\IEEEauthorblockA{Centre for Telecommunications Studies (CETUC)}\\
Pontifical Catholic University of Rio de Janeiro, Brazil \\
E-mail: \texttt{\{kaloorazi,delamare\}@cetuc.puc-rio.br}}

\markboth{Journal of \LaTeX\ Class Files,~Vol.~X, No.~X, January~XXXX}%
{Shell \MakeLowercase{\textit{et al.}}: Bare Demo of IEEEtran.cls for Journals}

\maketitle

\begin{abstract}
Low-rank matrix approximation plays an increasingly important role in signal and image processing applications. This paper presents a new rank-revealing decomposition method called randomized rank-revealing UZV decomposition (RRR-UZVD). RRR-UZVD is powered by randomization to approximate a low-rank input matrix. Given a large and dense matrix ${\bf A} \in \mathbb R^{m \times n}$ whose numerical rank is $k$, where $k$ is much smaller than $m$ and $n$, RRR-UZVD constructs an approximation $\hat{\bf A}$ such as $\hat{\bf A}={\bf UZV}^T$, where ${\bf U}$ and ${\bf V}$ have orthonormal columns, the leading-diagonal block of ${\bf Z}$ reveals the rank of $\bf A$, and its off-diagonal blocks have small $\ell_2$-norms. RRR-UZVD is simple, accurate, and only requires a few passes through $\bf A$ with an arithmetic cost of $O(mnk)$ floating-point operations. To demonstrate the effectiveness of the proposed method, we conduct experiments on synthetic data, as well as real data in applications of image reconstruction and robust principal component analysis.
\end{abstract}


\IEEEpeerreviewmaketitle

\section{Introduction}
\label{intro}
Low rank matrix approximation is approximating an input matrix by one of lower rank. The goal is to compactly represent the matrix with limited loss of information. Such a representation can provide a significant reduction in memory requirements as well as computational costs. Given a real large and dense $m \times n$ matrix ${\bf A}$ with numerical rank $k$ and $m \ge n$, its singular value decomposition (SVD) \cite{GolubVanLoan96} is written as follows: 
\begin{equation}
{\bf A} = {\bf U}_A{\bf \Sigma}_A{\bf V}_A^T= \begin{bmatrix} {{\bf U}_k \quad {\bf U}_0} \end{bmatrix}
  \begin{bmatrix}
       {\bf \Sigma}_k & 0  \\
       0 & {\bf \Sigma}_0
  \end{bmatrix}
  \begin{bmatrix}
       {{\bf V}_k \quad {\bf V}_0}
  \end{bmatrix}^T,
\label{eq1}
\end{equation}
where the matrices ${\bf U}_k \in \mathbb R^{m \times k}$, ${\bf U}_0 \in \mathbb R^{m \times n-k}$, ${\bf V}_k \in \mathbb R^{n \times k}$ and ${\bf V}_0 \in \mathbb R^{n \times n-k}$ are orthonormal, and ${\bf \Sigma}_k \in \mathbb R^{k \times k}$ and ${\bf \Sigma}_0 \in \mathbb R^{n-k \times n-k}$ are diagonal with entries $\rho_i$s as singular values. The SVD  constructs the best rank-$k$ approximation ${\bf B}$ to ${\bf A}$ \cite{GolubVanLoan96}, i.e., 
\begin{equation}
\begin{aligned}
&\ \|{\bf A} -  {\bf B}\|_2 = \rho_{k+1},  \\
&\ {\|{\bf A} - {\bf B}\|_F} = 
\Big(\rho_{k+1}^2 +...+ \rho_{n}^2 \Big)^{1/2},
\end{aligned}\label{equ9}
\end{equation} 
where  ${\bf B}= {\bf U}_k{\bf \Sigma}_k {\bf V}_k^T$, ${\|{\cdot}\|_2}$ and ${\|{\cdot}\|_F}$ denote the spectral norm and the Frobenius norm, respectively. In this paper we focus on the matrix $\bf A$ defined above. Another method to construct low-rank approximations of matrices is column-pivoted QR decomposition, or rank-revealing QR (RRQR) \cite{Chan87}. The RRQR \cite{Chan87} is a special QR with column pivoting (QRCP) which reveals the rank, i.e., the gap in the singular value spectrum, of the input matrix. The RRQR factors the matrix $\bf A$ as the product:
\begin{equation}
{\bf A} = {\bf Q}{\bf R}{\bf \Pi}^T= {\bf Q}
  \begin{bmatrix}
       {\bf R}_{11} & {\bf R}_{12}  \\
       {\bf 0} & {\bf R}_{22}
  \end{bmatrix}{\bf \Pi}^T,
\label{equTwolem1}
\end{equation}
where ${\bf Q} \in \mathbb R^{m \times n}$ is an orthonormal matrix, ${\bf R} \in \mathbb R^{n \times n}$ is an upper triangular matrix where ${\bf R}_{11} \in \mathbb R^{k \times k}$ is well-conditioned with $\rho_\text{min}({\bf R}_{11})= O(\rho_k)$, and ${\bf R}_{22} \in \mathbb R^{n-k \times n-k}$ has a sufficiently small $\ell_2$-norm. The matrix ${\bf \Pi} \in \mathbb R^{n \times n}$ is a permutation matrix. 

Low-rank matrices arise in many applications such as ranking and collaborative filtering \cite{Srebro2005}, background modeling \cite{WPMGR2009,BouPHV14,MFKDeTSP18,8425659,WWCPSH18}, image reconstruction \cite{MFKDeJSTSP18,DuerschGu2017,LiuSDWW8463507}, system identification \cite{FazelPST13}, Internet Protocol network anomaly detection \cite{KaDeICASSP17,NHussainP18}, sensor and multichannel signal processing \cite{DeSa2009}, and biometrics \cite{VictorBS02}. These two widely used methods, however, are computationally prohibitive for large matrices. In addition, their computations using standard schemes are challenging to be parallelized on advanced computational architectures \cite{HMT2009,Gu2015,AnztDGKLTY17}. Recently, low-rank approximation methods based on randomization have been developed 
\cite{FriezeKVS04,GoreinovTZ97,DrineasKM06,Sarlos06,Rokhlin09,HMT2009,BoutsidisG13,Gu2015,MFKDeTSP18}.
Randomized methods first transform an input matrix into a lower dimensional space by means of random matrices, and next apply traditional methods to further process the reduced-size matrix. These methods have been demonstrated to be remarkably efficient, accurate, and robust. The performances of randomized methods are known to be superior to those of the classical methods in many practical circumstances. The advantage of randomized methods over their classical counterparts is twofold: i) they are computationally more efficient, and ii) they readily lend themselves to a parallel implementation in order to exploit parallel architectures. 

\begin{itemize}
\item [] \textit{Structure}
\end{itemize}

The remainder of this paper is organized as follows. In Section \ref{Sec2ProbState} we discuss prior works, the problem this work is concerned with, and our contributions. In Section \ref{sec3RRR-UZVD} we describe our proposed method in detail. In Section \ref{sec4AnalysisUZVD} we present the mathematical analysis of RRR-UZVD. In Section \ref{sec6SimExp} we present and discuss our experimental results, and concluding remarks are given in Section \ref{sec7Conclu}.

\section{Prior Works and Problem Statement}
\label{Sec2ProbState}

Economical variants of SVD and RRQR are partial SVD and truncated RRQR \cite{GolubVanLoan96}. Both methods can compute an approximation of a matrix $\bf A$ with $O(mnk)$ floating-point operations (flops). Partial SVD is computed by invoking Krylov subspace methods, such as the Lanczos and Arnoldi algorithms. These methods, however, have two disadvantages: i) inherently, they are numerically unstable \cite{GolubVanLoan96}, and ii) they are difficult to parallelize \cite{HMT2009,Gu2015}, making them unsuitable to apply on advanced computational platforms. On the other hand, the major shortcoming in computing RRQR is pivoting strategy in which the columns of the large input matrix need to be permuted, i.e., RRQR is not \textit{pass-efficient} \cite{HMT2009}. To address this concern, recently several RRQR algorithms based on randomization have been proposed in order to carry out the factorization with minimum communication costs \cite{DemGGX15,DuerschGu2017,MartinssonHQRRP2017}. However, as we will show, RRQR does not provide highly accurate approximations.   

After emerging the works in \cite{FriezeKVS04,GoreinovTZ97}, many randomized algorithms have been proposed for computing low-rank approximations of matrices. The algorithms in \cite{DrineasKM06,DeshpandeV2006,RudelsonV07}, first sample columns of an input matrix according to a probability proportional to either their $\ell_2$-norm or leverage scores, leading to a compact represention of the matrix. The submatrix is then used for further computation using deterministic algorithms such as the SVD and pivoted QR \cite{GolubVanLoan96} to construct the final low-rank approximation. Rokhlin et al. \cite{Rokhlin09} proposed to first apply a random Gaussian embedding matrix to compress the input matrix in order to obtain an orthnormal basis for the range of the matrix. Next, the matrix is transformed into a lower-dimensional subspace by means of the basis. The low-rank approximation is then obtained through computations by means of an SVD on the reduced-size matrix. Sarl{\'{o}}s \cite{Sarlos06} proposed a method based on the Johnson-Lindenstrauss (JL) lemma. He showed that random linear combinations of rows render a good approximation to a low-rank matrix. The work in \cite{ClarWood2017}, built on Sarl{\'{o}}s's idea, computed a low-rank approximation using subspace embedding. Halko et al. \cite{HMT2009} developed an algorithm based on randomized sampling techniques to compute an approximate SVD of a matrix. Their method, \textit{randomized SVD}, for the matrix ${\bf A}$ and integer $0<k\le \ell< n$ is computed as described in Alg. \ref{Alg1}.

\begin{algorithm}
\caption{Randomized SVD (R-SVD) \cite{HMT2009}}
\renewcommand{\algorithmicrequire}{\textbf{Input:}}
\begin{algorithmic}[1]
\REQUIRE ~~ 
 Matrix $\ {\bf A} \in \mathbb R^{m \times n}$, integers $k>0$, $k \le \ell<n$.
\renewcommand{\algorithmicrequire}{\textbf{Output:}}
\REQUIRE ~~ A rank-$\ell$ approximation.
  \STATE Draw a random matrix ${\bf \Gamma} \in \mathbb R^{n \times \ell}$;
  \STATE Compute the matrix product ${\bf Y} = {\bf A}{\bf \Gamma}$;
  \STATE Compute a QR factorization ${\bf Y} = {\bf Q}{\bf R}$; 
  \STATE Compute the matrix product ${\bf B} = {\bf Q}^T\bf A$;
  \STATE Compute an SVD ${\bf B} = \bar{\bf U} \bar{\bf \Sigma}\bar{\bf V}^T$;
  \STATE Form the low-rank approximation $\hat{\bf A}_\text{RSVD} = ({\bf Q} \bar{\bf U})\bar{\bf \Sigma}\bar{\bf V}^T$.
\end{algorithmic}\label{Alg1}
\end{algorithm}

The R-SVD computes an SVD of an $n \times \ell$ matrix. 
For large matrices, specifically when $m \approx n$, this computation can be burdensome in terms of the communication cost, i.e., data movement either between different levels of a memory hierarchy or between processors \cite{DemmGHL12}. This makes the R-SVD unsuitable on modern computational environments. To address this issue, in this paper we develop a randomized algorithm that only utilizes the QR factorization to construct the approximation. The proposed algorithm, due to recently developed Communication-Avoiding QR (CAQR) algorithms \cite{DemmGHL12}, can be organized to exploit modern architectures, thereby being computed efficiently.

Gu \cite{Gu2015} slightly modified the R-SVD algorithm and applied it to improve subspace iteration methods. The work in \cite{MFKDeTSP18} proposed an algorithm termed subspace-orbit randomized SVD (SOR-SVD). SOR-SVD alternately projects the input matrix onto its column and row space via a random Gaussian matrix. The matrix is then transformed into a lower dimensional subspace, and a truncated SVD follows to construct an approximation. The work in \cite{MFKDeJSTSP18} proposed a rank-revealing algorithm based on randomized techniques termed compressed randomized UTV (CoR-UTV) decomposition. CoR-UTV first compresses the matrix through approximate orthonormal bases. Next, a QRCP is applied on the compressed matrix to give the final low-rank approximation.

\begin{itemize}
\item [] \textit{Our Contributions}
\end{itemize}
Driven by recent developments, this paper introduces a rank-revealing algorithm based on randomized sampling techniques termed randomized rank-revealing UZV decomposition (RRR-UZVD). RRR-UZVD constructs an approximation $\hat{\bf A}$ to the matrix ${\bf A}$ such as:
\begin{equation}
\hat{\bf A}={\bf UZV}^T,
\label{eq_contUZV}
\end{equation}  
where ${\bf U}$ and ${\bf V}$ are orthonormal matrices, the leading-diagonal block of ${\bf Z}$ reveals the rank of $\bf A$, and its off-diagonal blocks are sufficiently small in magnitude. RRR-UZVD requires a few passes through $\bf A$, and runs in $O(mnk)$ flops. The main operations of RRR-UZVD consist of matrix-matrix multiplications and QR factorizations. Recently, Communication-Avoiding QR algorithms \cite{DemmGHL12} have been developed. These algorithms can carry out a QR computation with optimal communication costs. Thus, RRR-UZVD can be optimized for peak machine performance on high performance computing devices. We establish a theoretical analysis for RRR-UZVD in which the rank-revealing property of the algorithm is proved (Theorem \ref{Thrm1}). Moreover, we apply RRR-UZVD to reconstruct a low-rank image. We further apply RRR-UZVD, as a surrogate to the expensive SVD, to solve the robust principal component analysis (robust PCA) problem \cite{WPMGR2009,CSPW2009,CLMW2009} in applications of background/foreground separation in video sequences as well as removing shadows and specularities from face images.

\section{Randomized Rank-Revealing UZV Decomposition}
\label{sec3RRR-UZVD}  

A rank-revealing decomposition is any decomposition in which the rank of a matrix is revealed \cite{Chan87,StewartURV92,StewartULV93}. The cr\`{e}me de la cr\`{e}me of rank-revealing decompositions is the SVD \cite{GolubVanLoan96}. Other deterministic rank-revealing methods include rank-revealing QR decomposition \cite{Chan87}, URV decomposition (URVD) \cite{StewartURV92}, and ULV decomposition (ULVD) \cite{StewartULV93}. The drawback of these deterministic methods, however, includes computational costs, both arithmetic and, more importantly, communication \cite{DemmGHL12,GunterVan05,DuerschGu2017}. 

The randomized rank-revealing UZV decomposition (RRR-UZVD) furnishes information on singular values and singular vectors of an input matrix through randomized schemes. Given ${\bf A}$, and integer $0<k\le \ell<n$, RRR-UZVD computes an approximation $\hat{\bf A}$ to ${\bf A}$ such that:
\begin{equation}
\hat{\bf A}={\bf U}{\bf Z}{\bf V}^T = {\bf U} 
  \begin{bmatrix}
       {\bf Z}_k & \bf G  \\
       \bf H & \bf E
  \end{bmatrix}{\bf V}^T,
\label{equ19}
\end{equation} 
where ${\bf U}$ and ${\bf V}$ are orthonormal matrices of size $m \times \ell$ and $n \times \ell$, respectively. The matrix ${\bf Z}_k$ is of order $k$ and nonsingular. The diagonals of ${\bf Z}_k$ are estimations to the first $k$ singular values of $\bf A$. The matrices ${\bf G}$, ${\bf H}$, and ${\bf E}$ are of size ${k \times \ell-k}$, ${\ell-k \times k}$, and $\ell-k \times \ell-k$, respectively all having small $\ell_2$-norms. We term the diagonal elements of ${\bf Z}$, $\bf Z$-$values$ of the matrix ${\bf A}$. The RRR-UZVD is rank-revealer in that the submatrix ${\bf Z}_k$ reveals the numerical rank $k$ of $\bf A$, and the other submatrices of ${\bf Z}$ have sufficiently small $\ell_2$-norms; see Theorem \ref{Thrm1}. The definition given for RRR-UZVD is analogous to those of rank-revealing algorithms in the literature 
\cite{Chan87,StewartURV92,StewartULV93,ChandIpsen94,MFKDeJSTSP18}.
To compute RRR-UZVD, the input matrix is first transformed into a lower dimensional space by means of random sampling techniques. Next, the entries of reduced-size matrix is manipulated. Lastly, the matrix is projected to the original space. For the matrix $\bf A$, RRR-UZVD is constructed by taking the following six steps:

\begin{enumerate}
\item Form a real $n \times \ell$ random test matrix 
${\bf \Theta}$.
\item Compute the $m \times \ell$ matrix product: 
\begin{equation}
{\bf F} = {\bf A}{\bf \Theta}.
\label{eqZ1}
\end{equation}
The matrix $\bf F$ is a projection onto a subspace spanned by $\bf A$'s columns.
\item Compute the $n \times \ell$ matrix product: 
\begin{equation}
{\bf T} = {\bf A}^T{\bf F}.
\label{eqZ2}
\end{equation}
The matrix $\bf T$ is a projection onto a subspace spanned by $\bf A$'s rows.
\item Form QR factorizations of ${\bf F}$ and ${\bf T}$: 
\begin{equation}
{\bf F} = {\bf U}{\bf R},  \quad \text{and} \quad {\bf T} = 
{\bf V}{\bf S},
\label{eqUV}
\end{equation}
where ${\bf U} \in \mathbb R^{m \times \ell}$ and ${\bf V} \in \mathbb R^{n \times \ell}$ are orthonormal matrices, which approximate the bases for $\mathcal{R} ({\bf A})$ and $\mathcal{R}({\bf A}^T)$, respectively. The notation $\mathcal{R}(\cdot)$ refers to the range of a matrix. The matrices ${\bf R}$, ${\bf S} \in \mathbb R^{\ell \times \ell}$ are upper triangular.  
\item Compute the ${\ell \times \ell}$ matrix product:
\begin{equation}
{\bf Z}={\bf U}^T{\bf A}{\bf V}.
\label{eqZ}
\end{equation}
The matrix ${\bf Z}$ is constructed by right and left  multiplications of $\bf A$ through the approximate bases.
\item Form the low-rank approximation of $\bf A$ by projecting the reduced matrix back to the original space: 
\begin{equation}
\hat{\bf A} = {\bf U}{\bf Z}{\bf V}^T.
\end{equation}
\end{enumerate}

The RRR-UZVD needs three passes over the data. But, it is possible to modify it in order for the algorithm to revisit the matrix ${\bf A}$ only once. For this purpose, the matrix $\bf Z$ \eqref{eqZ} is approximated through the available matrices: both sides of ${\bf Z}$ are right-multiplied by ${\bf V}^T {\bf \Theta}$, thus ${\bf Z}{\bf V}^T{\bf \Theta} = {\bf U}^T{\bf A}{\bf V}{\bf V}^T{\bf \Theta}$. Considering ${\bf A}\approx {\bf A}{\bf V}{\bf V}^T$, and ${\bf F} = {\bf A}{\bf \Theta}$, ${\bf Z}$ is approximated by
\begin{equation}
{\bf Z}_\text{approx} = {\bf U}^T{\bf F}({\bf V}^T{\bf \Theta})^\dagger,
\label{Zapp}
\end{equation} 
where $\dagger$ refers to the Moore–Penrose inverse. There exist two concerns associated with the basic version of RRR-UZVD: i) RRR-UZVD may provide poor approximate (leading) singular values and singular vectors, compared to those of the SVD, and ii) the columns of $\bf U$ and $\bf V$ may not be in a contributing order. Accordingly, the $\bf Z$-$values$ may not be sorted in a decreasing order. We propose two methods to cope with these issues:
\begin{itemize}
  \item [1.] \textit{Power iterations.} We incoroprate a few steps of the power method \cite{Rokhlin09,HMT2009}, which can substantially improve the performance of RRR-UZVD.
  \item [2.] \textit{Column permutation.} We implement a column pivoting technique as follows: first, sort the $\bf Z$-$values$ in a non-increasing order, delivering a permutation matrix ${\bf P}$. Second, right-multiply the matrices $\bf U$ and $\bf V$ by ${\bf P}$ such as ${\bf U}_s = {\bf U P}$, and ${\bf V}_s = \bf VP$.
\end{itemize}

The modified RRR-UZVD is presented in Alg. \ref{Alg3}. Note that in Alg. \ref{Alg3} (when the power iteration technique is utilized), a non-updated $\bf T$ must be used to form the approximation ${\bf Z}_\text{approx}$.

We also make use of column permutation technique, a shortcoming of QRCP when implemented on parallel architectures. However, there are fundamental differences between QRCP and RRR-UZVD in this regard. In RRR-UZVD:
\begin{itemize}
\item The column permutation strategy is implemented on matrices that are much smaller than the input matrix since $\ell \ll \text{min} \{m,n\}$.
\item The column permutation technique does not need access columns of the input matrix itself, i.e., it does not need passes over the data.
\end{itemize}

\begin{algorithm}
\caption{The RRR-UZVD Algorithm}
\renewcommand{\algorithmicrequire}{\textbf{Input:}}
\begin{algorithmic}[1]
\REQUIRE ~~ 
 Data matrix $\ {\bf A} \in \mathbb R^{m \times n}$,
integers $k>0$, $k \le \ell<n$, and power iteration factor $q$.
\renewcommand{\algorithmicrequire}{\textbf{Output:}}
\REQUIRE ~~ A low-rank approximation.
  \STATE Draw a random test matrix ${\bf T} \in \mathbb R^{n \times \ell}$;
  \FOR{$j=$ 1 to $q+1$}
   \STATE Compute the matrix product ${\bf F}={\bf A}{\bf T}$; \\
   \STATE Compute the matrix product ${\bf T}={\bf A}^T{\bf F}$;
  \ENDFOR \\
  \STATE Form QR factorizations ${\bf F} = {\bf U}{\bf R}$ 
  and ${\bf T} = {\bf V}{\bf S}$; 
  \STATE Compute ${\bf Z}_\text{approx} = {\bf U}^T{\bf F}({\bf V}^T{\bf}{\bf T})^\dagger$, (or ${\bf Z}={\bf U}^T{\bf A}{\bf V}$);
  \STATE Carry out the column pivoting technique, delivering 
  ${\bf U}_s$, ${\bf V}_s$, ${\bf Z}_\text {approx} = {\bf U}_s^T{\bf Z}_1 ({\bf V}_s^T{\bf}{\bf Z}_2)^\dagger $;
  \STATE Construct the low-rank approximation $\hat{\bf A}={\bf U}_s{\bf Z}_\text {approx}  {\bf V}_s^T$.
\end{algorithmic}\label{Alg3}
\end{algorithm}

The key difference between RRR-UZVD with SOR-SVD \cite{MFKDeTSP18} and CoR-UTV \cite{MFKDeJSTSP18} is that both SOR-SVD and CoR-UTV apply a deterministic algorithm, the former applies the SVD and latter applies column-pivoted QR, on the reduced matrix in their procedures to factor an input matrix. However, RRR-UZVD only uses column permutation and matrix-matrix multiplications on the small matrix. The deterministic methods employed can be challenging to parallelize when implemented on high performance computing devices, while matrix-matrix multiplications can readily be implemented in parallel, as will be explained later. In addition, SOR-SVD produces a rank-$k$ while RRR-UZVD constructs a rank-$\ell$ approximation.

\section{Analysis of RRR-UZV Decomposition}
\label{sec4AnalysisUZVD}

In this section, we show that RRR-UZVD has the rank-revealing property, and discuss its computational cost.
\subsection{Rank-Revealing Property}
The RRR-UZVD, with partitioned matrices, takes the form:  
\begin{equation}
\hat{\bf A} = {\bf U}{\bf Z}{\bf V}^T= \begin{bmatrix} {{\bf U}_1 \quad {\bf U}_2} \end{bmatrix}
  \begin{bmatrix}
       {\bf Z}_k & {\bf G}  \\
       {\bf H} & {\bf E}
  \end{bmatrix}\begin{bmatrix}{{\bf V}_1 \quad {\bf V}_2} \end{bmatrix}^T,
\label{equ24}
\end{equation}
where ${\bf U}_1$ and ${\bf U}_2$ are orthonormal matrices of size $m \times k$ and $m \times \ell-k$, respectively. 
${\bf Z}_k$ is of order $k$ and invertible whose diagonals being approximations of $\ell$ leading singular values of $\bf A$. The matrices ${\bf V}_1$ and ${\bf V}_2$ are orthonormal of size $n \times k$ and ${n \times \ell-k}$, respectively. We demonstrate that ${\bf Z}_k $ reveals the numerical rank of $\bf A$, and submatrices ${\bf G} \in \mathbb R^{k \times \ell-k}$, ${\bf H} \in \mathbb R^{\ell-k \times k}$, ${\bf E} \in \mathbb R^{\ell-k \times \ell-k}$ are small in magnitude. We show that if approximate bases for $\mathcal{R} ({\bf A})$ and $\mathcal{R}({\bf A}^T)$ are obtained through the procedure described, the middle matrix $\bf Z$ has the rank-revealing property. The theorem below sets forth the rank-revealing property of the RRR-UZVD algorithm. This result is new.

\begin{theorem}
Suppose that ${\bf A}$ is a real $m \times n$ matrix, where $m \ge n$, and its numerical rank is $k$. Further, suppose that the SVD of $\bf A$ is defined as in \eqref{eq1}, and its RRR-UZVD is defined as in \eqref{equ24}. Then,  
\begin{equation}
\rho_\text{min}({\bf Z}_k)= O(\rho_k),
\label{equThr1_1}
\end{equation}
\begin{equation}
\|[{\bf H} \quad {\bf E}]\|_2 = O(\rho_{k+1}),
\label{equThr1_2}
\end{equation}
\begin{equation}
\|[{\bf G}^T \quad {\bf E}^T]^T\|_2 = O(\rho_{k+1}).
\label{equThr1_3}
\end{equation}
\label{Thrm1}
\end{theorem}

\textit{Proof.} The proof is given in the Appendix. 

\subsection{Computational Cost}
\label{secComComplex}

The basic version of RRR-UZVD in order to construct an approximation to $\bf A$ has the following costs: 
\begin{itemize}
\item Forming a random matrix $\bf \Phi$, e.g., standard Gaussian, in Step 1 costs $O(n\ell)$.
\item Computing the matrix product in Step 2 costs $O(mn\ell)$.
\item Computing the matrix product in Step 3 costs $O(mn\ell)$.
\item Forming the QR factorizations in Step 4 costs $O(m\ell^2 + n\ell^2)$.
\item Computing the matrix product in Step 5 costs $O(mn\ell + m\ell^2)$ (if in this step the matrix $\bf Z$ is estimated by ${\bf Z}_\text{approx}$ \eqref{Zapp}, the cost would be $O(m\ell^2 + n\ell^2 +\ell^3)$). 
\end{itemize}

The dominant cost of Step 1 through Step 6 occurs when $\bf A$ and ${\bf A}^T$ are multiplied by the corresponding matrices. Hence
\begin{equation}
C_\text{UZV} = O(mn\ell).
\label{equCost1}
\end{equation}
The cost of the algorithm increases when the modifications (power iteration and column permutation techniques) are applied.  The column reordering technique costs $O(m\ell)$. RRR-UZVD needs either three or two passes (considering $\bf Z$ being approximated by ${\bf Z}_\text{approx}$) over data to decompose $\bf A$. With the power method used (Alg. \ref{Alg3}), RRR-UZVD needs either $(2q+3)$ passes over data with flop count of $(2q+3)C_\text{UZV}$, or $(2q+2)$ passes over data (considering $\bf Z$ being approximated by ${\bf Z}_\text{approx}$) with flop count of $(2q+2)C_\text{UZV}$. The sample size parameter $\ell$ is usually close to the exact numerical rank $k$. 

To compute any algorithm one needs to consider two costs: i) arithmetic, that is, the number of flops, and ii) communication, that is, data flow either between different levels of a memory hierarchy or between processors \cite{DemmGHL12}. The communication cost becomes considerably more expensive compared to the arithmetic when carrying out the computations for a large data matrix stored externally on modern computational platforms  \cite{DemmGHL12,Dongarra17}.
The RRR-UZVD carries out several matrix-matrix multiplications which can be readily implemented in parallel. RRR-UZVD further carries out QR factorizations on two $m\times \ell$ and $n\times \ell$ matrices. The R-SVD \cite{HMT2009}, however, carries out one QR factorization on a matrix of size $m\times \ell$ and one SVD on a matrix of size $n\times \ell$. Recently, Communication-Avoiding QR (CAQR) algorithms \cite{DemmGHL12} have been developed, which carry out the orthogonal triangularization with optimal communication costs. However, computing an SVD using  standard schemes are difficult to parallelize \cite{HMT2009,Gu2015,AnztDGKLTY17}. Hence, the operations of RRR-UZVD can be arranged to produce a low-rank approximation with optimal communication costs. This is an advantage of RRR-UZVD over R-SVD. Compared to CoR-UTV \cite{MFKDeJSTSP18}, RRR-UZVD only makes use of the matrix-vector multiplication to manipulate the reduced matrix $\bf Z$. While, CoR-UTV employs a deterministic QRCP for processing the compressed matrix, which may impose considerable communication cost on the algorithm. This is an advantage of RRR-UZVD over CoR-UTV.

\section{Simulations}
\label{sec6SimExp}

In this section, we assess the empirical performance of the RRR-UZVD algorithm. We first, through numerical examples, show that RRR-UZVD i) is a rank-revealing algorithm, and ii) furnishes estimates of singular values, which with remarkable fidelity track the exact singular values of the input matrix. To make a fair judgment of the behaviour of RRR-UZVD, we include QR with column pivoting (QRCP), the optimal SVD, R-SVD \cite{HMT2009}, and CoR-UTV \cite{MFKDeJSTSP18} in our comparison. Next, we treat an image reconstruction problem; we reconstruct a low-rank image of a differential gear through RRR-UZVD. Lastly, we devise a robust PCA algorithm by using RRR-UZVD. We experimentally investigate the efficacy and efficiency  of the proposed algorithm on synthetic and real-time data. The experiments were conducted in MATLAB on a PC with a 3 GHz intel Core i5-4430 processor and 8 GB of RAM. To compute RRR-UZVD and R-SVD, we generate a random matrix whose entries are independent, identically distributed Gaussian random variables of zero mean and unit variance.

\subsection{Rank-Revealing Property and Singular Value Approximation}
\label{subsecRMs}
We first assess the performance of RRR-UZVD on synthetically generated data. For the first test, we construct a rank-$k$ matrix ${\bf A}$ of order 1000 perturbed with Gaussian noise such that $\bf A ={\bf A}_1 +{\bf A}_2$ \cite{StewartQLP}. ${\bf A}_1={\bf U}_{{\bf A}_1}{\bf \Sigma}_{{\bf A}_1}{\bf  V}_{{\bf A}_1}^T$, where ${\bf U}_{{\bf A}_1}$ and ${\bf V}_{{\bf A}_1}$ are random orthonormal matrices, ${\bf \Sigma}_{{\bf A}_1}$ only has non-zero elements $\rho_i$s on the diagonal as singular values which decrease linearly from $1$ to $10^{-9}$, and $\rho_{k+1}=...=\rho_{1000}=0$. The matrix ${\bf A}_2$ is Gaussian normalized to have $\ell_{2}$-norm $\text{gap}\times \rho_k$. We set $k=20$, and gap = 0.15.

In the second test, we consider a challenging matrix called the devil's stairs \cite{StewartQLP}; a matrix ${\bf A}$ of order 1000 whose singular value spectrum has multiple gaps. The singular values are arranged analogous to a descending staircase. Each step comprises $d=10$ equal singular values. 

\begin{figure}[t]
\begin{center}
\input{GraphsUZV_2018/r_Stewart_2018}
\captionsetup{justification=centering,font=scriptsize}
\caption{Comparison of singular values of matrices considered computed by the SVD, QRCP, R-SVD \cite{HMT2009}, 
CoR-UTV \cite{MFKDeJSTSP18}, and proposed RRR-UZVD (Alg. \ref{Alg3}).} 
\label{fig:SteMat}       
\end{center}
\end{figure}
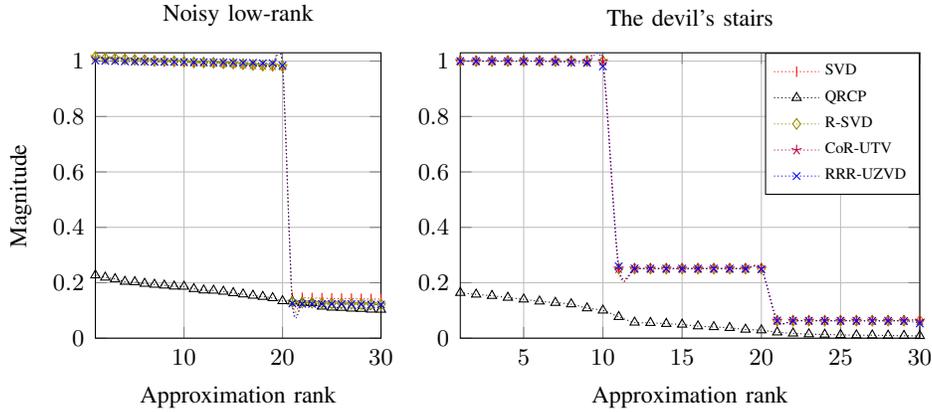

We factor the test matrices through RRR-UZVD (Alg. \ref{Alg3}), the SVD, QRCP, R-SVD, and CoR-UTV. We compare the results to assess the quality of singular values. For the randomized methods, we (arbitrarily) set $\ell=2k$, and the power method factor $q=1$. The results are shown in Fig. \ref{fig:SteMat}. We make two observations: i) for the first test matrix, RRR-UZVD strongly reveals the gap in $\rho_{20}$ and $\rho_{21}$, and provides estimations to singular values that show no loss of accuracy compared with the optimal SVD, while QRCP fails to reveal the rank. ii) For the second test matrix, RRR-UZVD reveals the gaps between the singular values and, moreover, excellently tracks them. RRR-UZVD furnishes highly accurate approximations to the singular values of the devil's stairs, whereas QRCP fails to reveal the gaps, provides poor estimations to the singular values, and is not able to track them. 

\subsection{Low-Rank Image Reconstruction}
\label{subImgRec}
In this experiment, we reconstruct a gray-scale image of a differential gear with dimension $1280\times 804$   \cite{DuerschGu2017} in order to evaluate the quality of  approximation constructed by RRR-UZVD.  In our comparison, we consider truncated QRCP, R-SVD, and the truncated SVD of PROPACK package \cite{Larsen98}. The PROPACK
function efficiently computes a specified number of leading singular pairs, suitable for approximating large matrices.
\begin{figure}
\begin{subfigure}{.5\textwidth}
  \includegraphics[width=0.9\textwidth,height=6.5cm]{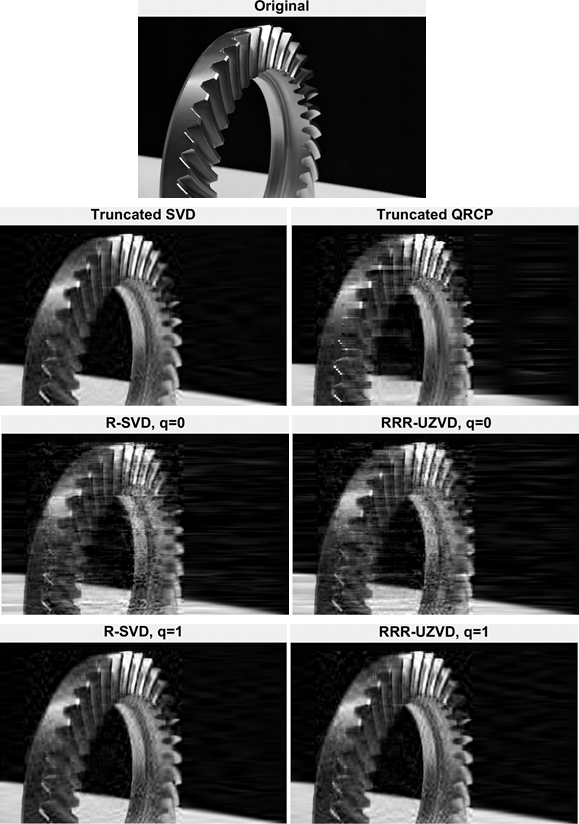}
  \caption{$Rank=25$ reconstruction}
\label{fig:Gear25}
\end{subfigure}%
\begin{subfigure}{.5\textwidth}
\includegraphics[width=0.9\textwidth,height=6.5cm]{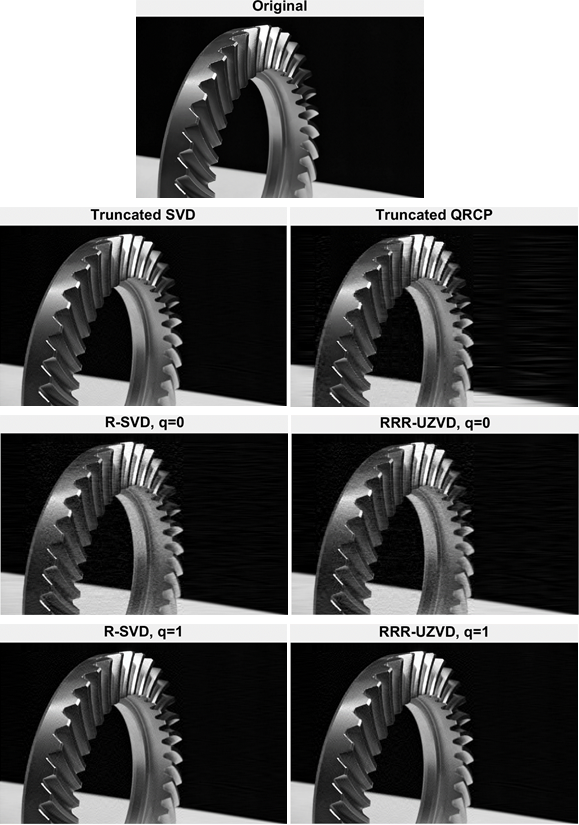}
  \caption{$Rank=85$ reconstruction}
\label{fig:Gear85}
\end{subfigure}
\caption{Low-rank image reconstruction through the truncated SVD \cite{Larsen98}, truncated QRCP, R-SVD \cite{HMT2009}, and RRR-UZVD with (a) $rank=25$, and (b) $rank=85$.}
\label{fig:ImgRec}
\end{figure}

\begin{figure}
\begin{subfigure}{.5\textwidth}
\input{GraphsUZV_2018/r_FroErr_Gear_UZV}
\captionsetup{justification=centering,font=scriptsize}
\caption{Approximation error} 
\label{fig:Fro_Gear}  
\end{subfigure}%
\begin{subfigure}{.44\textwidth}
\input{GraphsUZV_2018/r_Time_Gear_UZV}
\captionsetup{justification=centering,font=scriptsize}
\caption{Computational time} 
\label{fig:Time_Gear}
\end{subfigure}
\caption{(a) Approximation errors by the methods studied in reconstructing the low-rank differential gear image. (b) Runtime in seconds for different methods.}
\label{fig:FroTime}
\end{figure}
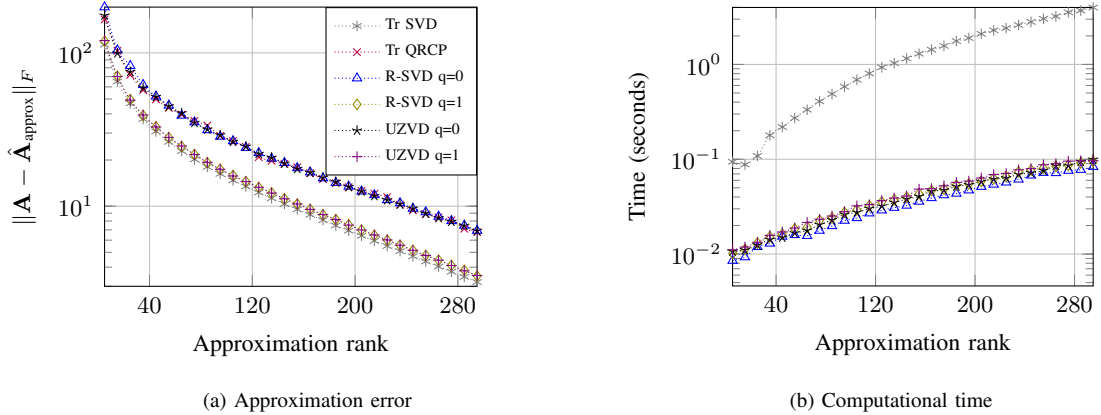

The results are displayed in Figs. \ref{fig:ImgRec} and \ref{fig:FroTime}. Fig. \ref{fig:ImgRec} shows the differential gear images reconstructed using $rank=25$ and $rank=85$, respectively via the methods considered. From Fig. \ref{fig:Gear25}, RRR-UZVD and R-SVD with no power iterations ($q=0$) demonstrate the poorest qualities. A better reconstruction is shown by the truncated QRCP. However, RRR-UZVD and R-SVD using one step of power iterations ($q=1$) construct approximations as accurate as the truncated SVD. In this case, these methods outperform the truncated QRCP. The rank-$85$ approximations in Fig. \ref{fig:Gear85} show tiny artifacts in images reconstructed through RRR-UZVD and R-SVD with $q=0$ as well as truncated QRCP. If one step of power iterations is incorporated, images reconstructed by truncated SVD, RRR-UZVD and R-SVD are visually identical. 

Fig. \ref{fig:Fro_Gear} shows the approximation error in terms of the Frobenius norm against the approximation rank. The error is calculated as follows:
\begin{equation}
\zeta = \|{\bf A} - \hat{\bf A}_{\text{approx}}\|_F,
\label{eq_ApprErr2}
\end{equation} 
where $\hat{\bf A}_{\text{approx}}$ is a low-rank approximation produced by each algorithm. In Fig. \ref{fig:Time_Gear} the execution times of the considered algorithms are compared, except for truncated QRCP; we have excluded this algorithm because there does not exist an optimized LAPACK function to compute QRCP with an assigned rank. The figure shows how the runtime of truncated SVD considerably grows as the approximation rank increases. While, one step of a power iteration together with column pivoting technique barely adds to the runtime of RRR-UZVD. This demonstrates the RRR-UZVD ability to provide comparable approximations with truncated SVD with much lower computational cost. R-SVD with $q=0$ appears to be slightly faster than RRR-UZVD with $q=0$ in our experiment. However, as RRR-UZVD's operations can be arranged to be carried out with minimum communication costs (subsection \ref{secComComplex}), we expect that on current and future advanced computational devices the RRR-UZVD algorithm to be faster than both truncated SVD and R-SVD, where the cost of communication is a serious bottleneck on performance of any algorithm.

\subsection{Robust PCA with RRR-UZVD}
\label{subrpca}
Robust PCA \cite{WPMGR2009,CSPW2009,CLMW2009} represents an $m \times n$ low-rank matrix ${\bf A}$ whose fraction of entries being grossly perturbed as a linear combination of a clean low-rank matrix ${\bf B}$ and a sparse matrix of errors ${\bf C}$ such as ${\bf A=B+C}$. Robust PCA solves the convex program:  
\begin{equation}
\begin{aligned}
&{\text{minimize}_{\bf(B, C)}} \ {\|{\bf B}\|_* + \gamma\|{\bf C}\|_1} \\
&{\text{subject to}} \ {\bf A} = {\bf B} + {\bf C},
\end{aligned}\label{equV1}
\end{equation} 
where ${\|\mbox{\bf M}\|_*}  \triangleq \sum_i\rho_i (\mbox{\bf M}) $ is defined as the nuclear norm of any matrix $\mbox{\bf M}$, ${\|\mbox{\bf M}\|_1} \triangleq \sum_{ij} |\mbox{\bf M}_{ij}|$ is defined as the $\ell_{1}$-norm of $\mbox{\bf M}$, and $\gamma>0$ is a weighting parameter. The method of augmented Lagrange multipliers (ALM) \cite{YY2009} is an efficient method to solve \eqref{equV1}, which yields the optimal solution. However, it has a major bottleneck which is computing an SVD at each iteration for approximating the low-rank component $\bf B$ of $\bf A$. The work in \cite{LLS2011} proposes several techniques such as predicting the dimension of principal singular space, a continuation method \cite{Toh2010}, and a truncated SVD by making use of PROPACK package \cite{Larsen98} to speed up the convergence of the ALM method. The bottleneck of the proposed algorithm \cite{LLS2011}, however, is that the employed truncated SVD \cite{Larsen98} uses the lanczos method that is i) numerically unstable, and ii) has poor performance on modern computing devices because of limited data reuse in its operations \cite{GolubVanLoan96,HMT2009,Gu2015}. To address this issue, by retaining the original objective function \cite{WPMGR2009,CSPW2009,CLMW2009,LLS2011}, we apply RRR-UZVD to solve the robust PCA problem. We further adopt the continuation technique \cite{Toh2010,LLS2011}. We call our proposed method \texttt{ALM-UZVD} given in Alg. \ref{Alg_ALM-UZVD}. In Alg. \ref{Alg_ALM-UZVD}, for any matrix $\bf M$ with an RRR-UZVD described in Section \ref{sec3RRR-UZVD}, $\mathcal{Z}_\nu ({\bf M})$ denotes a UZV thresholding operator defined as:  
\begin{equation}
\mathcal{Z}_\nu({\bf M})={\bf U}(:,1:s){\bf Z}(1:s,:){\bf V}^T,
\end{equation} 
where $s$ is the number of diagonal elements of $\bf Z$ larger than $\nu$, $\mathcal{S}_\nu (x) = {\text{sgn}(x)\text{max}}(|x| - \nu, 0)$ is a shrinkage operator, and $\gamma$, $\eta_0$, ${\bar \eta}$, $\tau$, ${\bf Y}_0$, and ${\bf C}_0$ are initial values. 

\begin{algorithm}
\caption{Robust PCA solved via \texttt{ALM-UZVD}}
\renewcommand{\algorithmicrequire}{\textbf{Input:}}
\begin{algorithmic}[1]
\REQUIRE ~~ 
Matrix ${\bf A}, \gamma, \eta_0,{\bar \eta}, \tau, {\bf Y}_0, {\bf C}_0, i=0$;
 \renewcommand{\algorithmicrequire}{\textbf{Output:}} 
 \REQUIRE ~~ 
 Low-rank plus sparse matrix
\WHILE {the algorithm does not converge}
        \STATE Compute ${\bf B}_{i+1}= \mathcal{Z}_{\eta_i^{-1}}
        ({\bf A} - {\bf C}_i +\eta_i^{-1} {\bf Y}_i)$;
        \STATE Compute ${\bf C}_{i+1} = \mathcal{S}_{\gamma\eta_i^{-1}}({\bf A} - {\bf B}_{i+1} +\eta_i^{-1} {\bf Y})$;
        \STATE Compute${\bf Y}_{i+1}={\bf Y}_i +\eta_i({\bf A} - {\bf B}_{i+1} - {\bf C}_{i+1})$;
        \STATE Update $\eta_{i+1} = \text{max}(\tau\eta_i, {\bar \eta})$;
\ENDWHILE
\RETURN $\bf B^*$ and $\bf C^*$ 
\end{algorithmic}\label{Alg_ALM-UZVD}
\end{algorithm}

\subsubsection{Synthetic Low-Rank and Sparse Matrix Recovery}
\label{subsecSyDataRec}

We construct a matrix $\bf A = B + C$ as a sum of a low-rank matrix ${\bf B} \in \mathbb R^{n \times n}$ and a sparse matrix of errors ${\bf C}\in \mathbb R^{n \times n}$. The matrix ${\bf B}$ is generated as a product of two standard Gaussian matrices ${\bf W}$, ${\bf Q} \in \mathbb R^{n \times k}$ such as ${\bf B}={\bf W}{\bf Q}^T$. The matrix ${\bf C}$ has $c$ non-zero entries drawn independently from the set $\lbrace$-100, 100$\rbrace$. We consider the rank $k = \text{rank}({\bf B})=0.05\times n$, and $c = \|{\bf C}\|_0=0.05\times n^2$, where $\|\cdot\|_0$ refers to the $\ell_0$-norm.

We apply the \texttt{ALM-UZVD} and efficiently implemented robust PCA algorithm in \cite{LLS2011}, hereafter \texttt{InexactALM}, to $\bf A$ to recover ${\bf B}^*$ and ${\bf C}^*$. The numerical results are summarized in 
Table \ref{TableFour}. We adopt the initial values proposed in \cite{LLS2011} in our experiments. The algorithms are terminated when ${\|{\bf A}-{\bf B}^{sol}-{\bf C}^{sol}\|_F}< 10^{-4}{{\|{\bf A}\|_F}}$ is satisfied, where $({\bf B}^{sol}, {\bf C}^{sol})$ is the solution pair of either algorithm. In Table \ref{TableFour}, $Time$ denotes the computational time in seconds, $Iter.$ denotes the number of iterations, and $\xi ={\|{\bf A}-{\bf B}^{sol}-{\bf C}^{sol}\|_F}/{{\|{\bf A}\|_F}}$ is the relative error. RRR-UZVD needs a preassigned rank $\ell$ to carry out the decomposition. We arbitrarily set $\ell=2k$, and also $q=2$. 

Judging from the results presented in Table \ref{TableFour}, we make the following observations: i) in all cases, \texttt{ALM-UZVD} successfully detects the exact rank $k$ of the matrices , ii) \texttt{ALM-UZVD} provides the optimal solution, while it needs one more iteration in comparison with \texttt{InexactALM}, and iii) \texttt{ALM-UZVD} outperforms \texttt{InexactALM} in runtime, with speedups of up to 5 times.

\begin{table}[!htb]
\centering
\caption{Numerical results of \texttt{InexactALM} \cite{LLS2011} and proposed \texttt{ALM-UZVD} (Alg. \ref{Alg_ALM-UZVD}) for synthetic matrix recovery.}
\begin{tabular}
{p{0.3cm} p{0.3cm} p{0.4cm} p{1.7cm} p{0.2cm} p{0.15cm} p{0.4cm} p{0.4cm} p{0.1cm}}
\noindent\rule{7.9cm}{0.4pt}\\
$n$ & $k$ & $c$ & Algorithms & 
$\hat{k}$ & $\hat{c}$ & Time & Iter.& $\xi$ \\
\noindent\rule{7.9cm}{0.4pt}\\
1e3& 50 & 5e4 & 
\begin{tabular}{|c p{0.2cm} p{0.3cm} p{0.3cm} p{0.1cm} p{0.9cm}}
\texttt{ALM-UZVD} & 50& 5e4 & 0.6 & 10 & 4.2e-5 \\
\texttt{InexactALM} & 50& 5e4 & 2.5 & 9 & 3.1e-5  \\
\end{tabular} \\
\\
2e3& 100 & 2e5 &
\begin{tabular}{|c p{0.2cm} p{0.3cm} p{0.3cm} p{0.1cm} p{0.9cm}}
\texttt{ALM-UZVD} & 100& 2e5 & 4.4 & 10 & 4.1e-5 \\
\texttt{InexactALM} & 100& 2e5 & 17.6 & 9 & 4.9e-5  \\
\end{tabular} \\
\\
3e3& 150 & 45e4 &
\begin{tabular}{|c p{0.2cm} p{0.3cm} p{0.3cm} p{0.1cm} p{0.9cm}}
\texttt{ALM-UZVD} & 150& 45e4 & 10.5 & 10 & 5.3e-5   \\
\texttt{InexactALM} & 150& 45e4 & 52.9 & 9 & 5.2e-5  \\
\end{tabular} \\
\vspace{-.15cm}
\noindent\rule{7.9cm}{0.4pt} 
\end{tabular}
\label{TableFour}
\end{table}

\subsubsection{Background Subtraction in Surveillance Video}
\label{subsecBSub}
In this experiment, we employ \texttt{ALM-UZVD} in order to separate the background and foreground of a video sequence. We use one surveillance video from \cite{LHGT2004}. The video has 200 grayscale frames each with dimension ${256 \times 320}$, taken in a shopping mall. We form an $81920 \times 200$ matrix ${\bf A}$ through stacking individual frames as its columns.
 
\begin{figure}
\centering
\begin{subfigure}{.4\textwidth}
  \centering
  \includegraphics[width=0.85\textwidth,height=4.5cm]{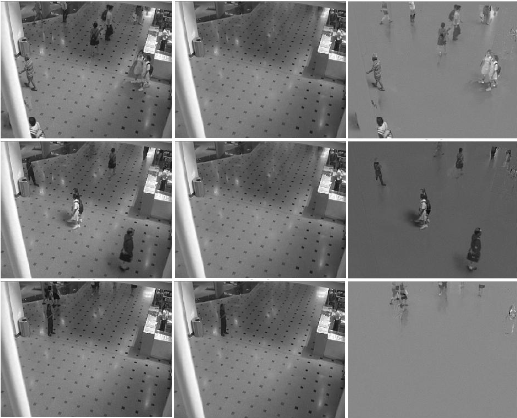}
  \caption{Shopping mall}
  \label{fig:sub1}
\end{subfigure}%
\begin{subfigure}{.4\textwidth}
  \centering
  \includegraphics[width=0.83\textwidth,height=4.5cm]{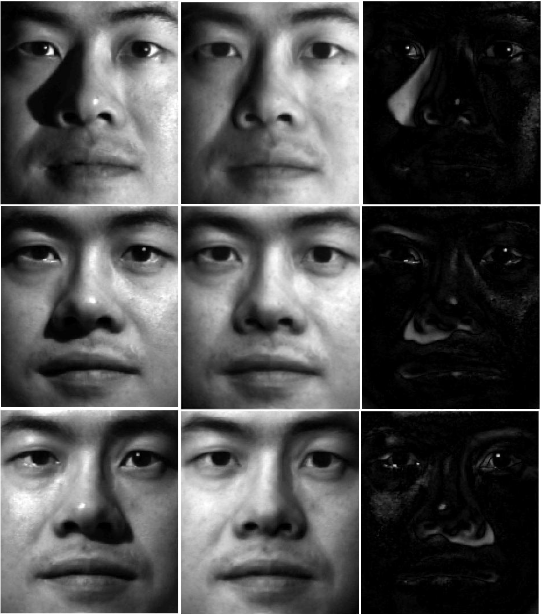}
  \caption{Yale B02 face}
  \label{fig:sub2}
\end{subfigure}
\caption{(a) Original video frames are shown in column 1. Recovered backgrounds ${\bf B}^*$ and foregrounds ${\bf C}^*$, by \texttt{ALM-UZVD} are shown in columns 2 and 3, respectively. (b) Cropped images of a face under different illuminations are shown in column 1. Clean images and errors associated with shadows and specularities recovered by \texttt{ALM-UZVD} are shown in columns 2 and 3, respectively.}
\label{fig:BackModel}
\end{figure}

To approximate the low-rank component (background), the RRR-UZVD incorporated in \texttt{ALM-UZVD} requires a preassigned rank $\ell$. We use the following bound \cite{GolubVanLoan96} that relates the rank $k$ of any matrix $\bf M$ with the Frobenius and nuclear norms in order to determine $\ell$: 
\begin{equation}
\sqrt{k} \ge \dfrac{\|{\bf M}\|_*}{\|{\bf M}\|_F}.
\label{equRank}
\end{equation}
We assign $\ell=k + p$, where $k$ is the minimum value satisfying \eqref{equRank}, and $p = 2$ is an oversampling factor. We again set $q=2$ for RRR-UZVD. Some video frames with recovered backgrounds and foregrounds are shown in Fig. \ref{fig:sub1}. It is observed that \texttt{ALM-UZVD} recovers the low-rank and sparse components of the video successfully. Table \ref{TableSix} presents the results, showing that \texttt{ALM-UZVD} outperforms \texttt{InexactALM} in terms of runtime.

\subsubsection{Removing Shadows and Specularities from Face Images}
\label{subsecShRemoval}
In this experiment, we apply \texttt{ALM-UZVD} in order to remove shadows and specularities of face Images. We use a face image of dimenstion ${192 \times 168}$ with a total of 64 illuminations from the Yale B face database \cite{Georghiades2001}. The individual images are stacked as the columns to form a $32256 \times 64$ matrix ${\bf A}$. We apply \texttt{ALM-UZVD} to $\bf A$. The results are shown in Fig. \ref{fig:sub2}, and Table \ref{TableSix}. It is observed that the sparse components contain the shadows and specularities of the face images that have been effectively extracted by \texttt{ALM-UZVD}. The \texttt{ALM-UZVD} method, we conclude, successfully recovers the face images taken under distant illumination from the dataset considered nearly two times faster than \texttt{InexactALM}. 

\begin{table}[!htb]
\centering
\caption{Numerical results of \texttt{InexactALM} \cite{LLS2011} and \texttt{ALM-UZVD} (Alg. \ref{Alg_ALM-UZVD}) for real data recovery.}
\begin{tabular}{p{2cm} p{2.5cm} p{2.5cm}} 
\noindent\rule{8cm}{0.4pt}\\
Dataset  
 & \begin{tabular}{p{0.4cm} p{0.5cm} p{0.5cm}}
   \multicolumn{3}{c}{\texttt{ALM-UZVD}} \\
   \hline
   Time & Iter. & $\xi$ \\
   \end{tabular}
   &
   \begin{tabular}{p{0.42cm} p{0.5cm} p{0.5cm}}
 \multicolumn{3}{c}{\texttt{InexactALM}} \\
\hline
 Time & Iter. & $\xi$ \\
 \end{tabular}\\
\noindent\rule{8cm}{0.4pt}\\
Shopping mall  & 
\begin{tabular}{p{0.4cm} ccp{0.75cm}c c}
15.8 &23 &8.1e-5 & $\mspace{2mu}$ 47.8 & $\mspace{-18mu}$ 23 &7.1e-5\\
\end{tabular} \\\\
Yale B02 & 
\begin{tabular}{cccccc}
2.0 & 21 &9.1e-5& $\mspace{5mu}$ 3.5 & 21 & $\mspace{-6mu}$ 7.5e-5\\
\end{tabular} \\
\vspace{-0.2cm}
\noindent\rule{8cm}{0.4pt}\\
\end{tabular}
\label{TableSix}
\end{table}

\section{Conclusion}
\label{sec7Conclu} 

In this paper we proposed RRR-UZVD. The RRR-UZVD method makes use of randomized sampling to construct an approximation to a low-rank matrix. We provided theoretical analysis for RRR-UZVD, showing the algorithm is rank-revealer. Through numerical experiments we showed that RRR-UZVD considerably outperforms QRCP in revealing the numerical rank of a low-rank matrix, and provides results with no loss of accuracy compared with those of the optimal SVD. The performance of RRR-UZVD exceeds that of QRCP in computing low-rank approximation when the power iteration and column permutation schemes are employed. In this case, RRR-UZVD furnishes results as accurate as those of the SVD. Compared to the SVD and QRCP, RRR-UZVD is more efficient in terms of arithmetic cost. In addition, RRR-UZVD can take advantage of modern computer architectures better than SVD, QRCP, as well as competing R-SVD and CoR-UTV. We further applied RRR-UZVD in applications of robust PCA. Our results show that \texttt{ALM-UZVD} renders the optimal solution, and is considerably faster than the efficiently implemented \texttt{InexactALM} method.

\section{Appendix}
Proof of Theorem \ref{Thrm1}.

To prove the first relation in \eqref{equThr1_1}, consider ${\bf B}\in \mathbb R^{m \times m}$ and ${\bf C} \in \mathbb R^{n \times n}$ to be orthonormal matrices (i.e., ${\bf B}^T{\bf B}={\bf I}_m$ and ${\bf C}^T{\bf C}={\bf I}_n$, where ${\bf I}_s$ denotes an $s \times s$ identity matrix). The matrices ${\bf A}$ and $\bf BAC$ have exactly the same singular values \cite{GolubVanLoan96}. Let a matrix ${\bf D}$ be a submatrix of $\bf BAC$ such that:
\begin{equation}
{\bf D} = {\bf B}_1{\bf A}{\bf C}_1,
\label{equ1Th}
\end{equation} 
where ${\bf B}_1\in \mathbb R^{k \times m}$ consists of the first $k$ rows of $\bf B$, and ${\bf C}_1\in \mathbb R^{n \times k}$ consists of the first $k$ columns of $\bf C$. Suppose that ${\bf D}$ has singular values such as $\delta_1 \ge \delta_2 \ge ...\ge \delta_k$. We form the following product: 
\begin{equation}
{\bf DD}^T = {\bf B}_1{\bf MM}^T{\bf B}_1^T,
\end{equation} 
where the matrix ${\bf M} = {\bf A}{\bf C}_1$ is of size $m\times k$. Thus, ${\bf DD}^T$ is a $k\times k$ principal submatrix of the $m\times m$ symmetric ${\bf B}{\bf MM}^T{\bf B}^T$. Suppose that ${\bf M}$ has singular values such as $\mu_1 \ge \mu_2 \ge ...\ge \mu_k$. Thus, the eigenvalues of ${\bf MM}^T$ satisfy
\begin{equation}
\mu_1^2 \ge \mu_2^2 \ge ...\ge \mu_k^2 \ge \mu_{k+1}^2 = ... = \mu_m^2=0.
\end{equation}
By invoking the formulas that relate the eigenvalues of a principal submatrix with those of a symmetric matrix  \cite{GolubVanLoan96}, we obtain
\begin{equation}
\mu_1^2 \ge \delta_1^2 \ge ...\ge \mu_k^2 \ge \delta_k^2.
\end{equation}
We now form the following product: 
\begin{equation}
{\bf M}^T{\bf M} = {\bf C}_1^T{\bf A}^T{\bf A}{\bf C}_1.
\label{eqCACA}
\end{equation} 
The nonzero eigenvalues of ${\bf M}^T{\bf M}$ coincide with those of ${\bf MM}^T$. By \eqref{eqCACA} we observe that ${\bf M}^T\bf M$ is a $k\times k$ principal submatrix of the $n\times n$ symmetric ${\bf C}^T{\bf A}^T{\bf A}{\bf C}$. Hence  
\begin{equation}
\rho_1^2 \ge \mu_1^2 \ge ...\ge \rho_k^2 \ge \mu_k^2,
\end{equation} 
and accordingly, 
\begin{equation}
\rho_1^2 \ge \mu_1^2 \ge \delta_1^2 \ge ...\ge \rho_k^2 \ge \mu_k^2 \ge\delta_k^2.
\end{equation}
Substituting ${\bf B}_1$ and ${\bf C}_1$ in \eqref{equ1Th} with ${\bf U}_1^T$ and ${\bf V}_1$ in \eqref{equ24}, respectively, results in 
\begin{equation}
\rho_k({\bf Z}) \le \rho_k.
\end{equation}

Note that for any decomposition to be rank-revealing, there must be a well-defined gap in the singular value spectrum of the input matrix \cite{Chan87,ChandIpsen94}.
 
To prove \eqref{equThr1_2}, with an SVD of $\bf A$ defined in \eqref{eq1}, let the matrix $\hat{\bf A}$ formed by RRR-UZVD have an SVD such that $\hat{\bf A} = \hat{\bf U}\hat{\bf \Sigma}\hat{\bf V}^T$. We write $\hat{\bf A} = \hat{\bf U}_k\hat{\bf \Sigma}_k\hat{\bf V}_k^T+ \hat{\bf U}_0 \hat{\bf \Sigma}_0\hat{\bf V}_0^T$. By the argument given in the proof of the first bound \eqref{equThr1_1}, the following relation must hold 
\begin{equation}
\rho_{k+1} = \|{\bf A}^T{\bf U}_0\|_2 \ge \hat{\rho}_{k+1} = 
\|\hat{\bf A}^T\hat{\bf U}_0\|_2.
\label{equ39}
\end{equation}  
We further have 
\begin{equation}
\hat{\rho}_{k+1} = 
\|\hat{\bf A}^T\hat{\bf U}_0\|_2 \le 
\|\hat{\bf A}^T{\bf U}_2\|_2.
\label{equ39_1}
\end{equation}  
The relation in \eqref{equ39_1} holds since ${\bf U}_2$ provides an approximation to $\hat{\bf U}_0$. However, in practice, combining the power iteration and column pivoting techniques leads to a ${\bf U}_2$ which spans the null space of $\hat{\bf A}^T$. Thus
\begin{equation}
\hat{\rho}_{k+1} \approx 
\|\hat{\bf A}^T{\bf U}_2\|_2.
\label{equ39_2}
\end{equation}  
By substituting $\hat{\bf A}$ \eqref{equ24} into \eqref{equ39_2}, it follows that
\begin{equation}
\begin{aligned}
 \hat{\rho}_{k+1} & { \approx \|\begin{bmatrix} {{\bf V}_1\quad{\bf V}_2} \end{bmatrix}
  \begin{bmatrix}
       {\bf Z}_k^T & {\bf H}^T  \\
       {\bf G}^T & {\bf E}^T
  \end{bmatrix}
  \begin{bmatrix}
       {\bf U}_1^T \\
       {\bf U}_2^T
       \end{bmatrix}}{\bf U}_2\|_2 \\            
& { \approx \|[{\bf H} \quad {\bf E}]\|_2}.
\end{aligned} 
\end{equation} 

To prove \eqref{equThr1_3}, with an analogous argument as in the proof of \eqref{equThr1_2}, we have
\begin{equation} 
\begin{aligned}
\rho_{k+1}  = \|{\bf A}{\bf V}_0\|_2 \ge \hat{\rho}_{k+1} = 
\|\hat{\bf A}\hat{\bf V}_0\|_2.
\end{aligned}\label{equ41*}
\end{equation}  
and consequently,
\begin{equation}
\hat{\rho}_{k+1} = 
\|\hat{\bf A}\hat{\bf V}_0\|_2 \le 
\|\hat{\bf A}{\bf V}_2\|_2.
\label{equ41}
\end{equation}
The relation in \eqref{equ41} holds since ${\bf V}_2$ is an approximation to $\hat{\bf V}_0$. In practice, combining the power iteration and column permutation techniques leads to a ${\bf V}_2$ which spans the null space of $\hat{\bf A}$. Thus
\begin{equation}
\begin{aligned}
\hat{\rho}_{k+1}  \approx
\|\hat{\bf A}{\bf V}_2\|_2.
\end{aligned}\label{equ41_1}
\end{equation}
By substituting $\hat{\bf A}$ \eqref{equ24} into \eqref{equ41_1}, 
it follows that 
\begin{equation}
\begin{aligned}
 \hat{\rho}_{k+1} & { \approx \|\begin{bmatrix} {{\bf U}_1\quad{\bf U}_2} \end{bmatrix}
  \begin{bmatrix}
       {\bf Z}_k & {\bf G}  \\
       {\bf H} & {\bf E}
  \end{bmatrix}
  \begin{bmatrix}
       {\bf V}_1^T \\
       {\bf V}_2^T
       \end{bmatrix}}{\bf V}_2\|_2 \\   
&       {\approx \|[{\bf G}^T \quad {\bf E}^T]^T\|_2}.
\end{aligned} 
\end{equation}
This completes the proof.




The authors would like to thank...

\ifCLASSOPTIONcaptionsoff
  \newpage
\fi

\bibliographystyle{IEEEtran}
\bibliography{mybibfileUZV}

\end{document}

%% file: GraphsUZV_2018/r_Stewart_2018.tex
%
%
%
\usetikzlibrary{positioning,calc}

\definecolor{mycolor1}{rgb}{0.00000,1.00000,1.00000}%
\definecolor{mycolor2}{rgb}{1.00000,0.00000,1.00000}%

\pgfplotsset{every axis label/.append style={font=\footnotesize},
every tick label/.append style={font=\footnotesize}
}

\begin{tikzpicture}[font=\footnotesize] 

\begin{axis}[%
name=ber,
width  = 0.23\columnwidth,
height = 0.23\columnwidth,
scale only axis,
xmin  = 1,
xmax  = 30,
xlabel= {Approximation rank},
xmajorgrids,
ymin=0,
ymax=1.03,
ylabel={Magnitude},
ymajorgrids,
title = {Noisy low-rank}
]

\addplot+[smooth,color=red,densely dotted, every mark/.append style={solid}, mark=|]
table[row sep=crcr]{
1	1.01484875150844\\
2	1.00995862938108\\
3	1.00927806416332\\
4	1.00695448462165\\
5	1.00483434196015\\
6	1.00357542841725\\
7	1.00054283264108\\
8	0.999933747229965\\
9	0.998889073445844\\
10	0.997383700740246\\
11	0.995262568075296\\
12	0.994926285874388\\
13	0.993982035422763\\
14	0.992212093777547\\
15	0.991278513188547\\
16	0.988568552281987\\
17	0.986155087500802\\
18	0.983951807268206\\
19	0.983327262961779\\
20	0.981665069421062\\
21	0.145867488994028\\
22	0.145612188285977\\
23	0.144753177420806\\
24	0.144236782930346\\
25	0.143189783014996\\
26	0.142820683010049\\
27	0.142493230054508\\
28	0.141699894195107\\
29	0.141301549921166\\
30	0.140981840676501\\
};

\addplot+[smooth,color=black,densely dotted, every mark/.append style={solid}, mark=triangle]
table[row sep=crcr]{
  1	0.227332213024648\\
2	0.219251369981432\\
3	0.213683857911715\\
4	0.204199916750572\\
5	0.202671420525675\\
6	0.196941326394703\\
7	0.193576705903990\\
8	0.190559098743500\\
9	0.186969282936563\\
10	0.186455949895567\\
11	0.177975303355588\\
12	0.172877921734879\\
13	0.172152583215336\\
14	0.168424837154943\\
15	0.163150964453923\\
16	0.158547855690776\\
17	0.153952448666359\\
18	0.149772975507999\\
19	0.145043313763871\\
20	0.134207562709647\\
21	0.133141058247291\\
22	0.128203659017668\\
23	0.127059717716467\\
24	0.114595748240451\\
25	0.110935558387716\\
26	0.110103874400236\\
27	0.107368711579790\\
28	0.106457799871307\\
29	0.103743025544514\\
30	0.103427151374411\\
  };

\addplot+[smooth,color=olive,densely dotted, every mark/.append style={solid}, mark=diamond]
  table[row sep=crcr]{
  1	1.01482777836683\\
2	1.00992958922637\\
3	1.00925698748761\\
4	1.00693104681132\\
5	1.00481606317181\\
6	1.00355771095704\\
7	1.00052720538931\\
8	0.999925401735529\\
9	0.998876952253299\\
10	0.997370578868139\\
11	0.995248176345529\\
12	0.994906556334548\\
13	0.993963986485288\\
14	0.992193031808204\\
15	0.991255343512511\\
16	0.988557010067719\\
17	0.986130748958646\\
18	0.983938779667039\\
19	0.983309463190147\\
20	0.981649137980620\\
21	0.130711918971353\\
22	0.128573264882789\\
23	0.127888176607429\\
24	0.126346203026367\\
25	0.125825024645770\\
26	0.125386723563492\\
27	0.123929232714362\\
28	0.123714789582857\\
29	0.123229053992272\\
30	0.121374713751041\\
  };
  
  \addplot+[smooth,color=olive,densely dotted, every mark/.append style={solid}, mark=star]
  table[row sep=crcr]{
  1	1.01482777836683\\
2	1.00992958922637\\
3	1.00925698748761\\
4	1.00693104681132\\
5	1.00481606317181\\
6	1.00355771095704\\
7	1.00052720538931\\
8	0.999925401735529\\
9	0.998876952253299\\
10	0.997370578868139\\
11	0.995248176345529\\
12	0.994906556334548\\
13	0.993963986485288\\
14	0.992193031808204\\
15	0.991255343512511\\
16	0.988557010067719\\
17	0.986130748958646\\
18	0.983938779667039\\
19	0.983309463190147\\
20	0.981649137980620\\
21	0.130711918971353\\
22	0.128573264882789\\
23	0.127888176607429\\
24	0.126346203026367\\
25	0.125825024645770\\
26	0.125386723563492\\
27	0.123929232714362\\
28	0.123714789582857\\
29	0.123229053992272\\
30	0.121374713751041\\
  };
  
\addplot+[smooth,color=blue,densely dotted, every mark/.append style={solid}, mark=x]
  table[row sep=crcr]{
1	1.00173334879797\\
2	1.00050345279511\\
3	0.999724359379005\\
4	0.998713887665144\\
5	0.998027909375524\\
6	0.997969633286014\\
7	0.997398666052992\\
8	0.997358123336721\\
9	0.997229612380842\\
10	0.996602147930927\\
11	0.995925973387042\\
12	0.995878036025569\\
13	0.995258498112571\\
14	0.995136362377284\\
15	0.994432726898652\\
16	0.993683577890047\\
17	0.993682077990525\\
18	0.992370113005767\\
19	0.992336558263670\\
20	0.984092352122066\\
21	0.125324457822016\\
22	0.123818312021468\\
23	0.123519884931627\\
24	0.123470797234303\\
25	0.123455729447670\\
26	0.123073957016453\\
27	0.122883597292440\\
28	0.122346701871804\\
29	0.121407726479048\\
30	0.121222766501217\\
  };  
\end{axis}

\begin{axis}[%
name=SumRate,
at={($(ber.east)+(30,0em)$)}, anchor= west,
width  = 0.37\columnwidth,
height = 0.23\columnwidth,
scale only axis,
xmin  = 1,
xmax  = 30,
xlabel= {Approximation rank},
xmajorgrids,
ymin=0,
ymax=1.03,
ylabel={},
ymajorgrids,
legend entries={SVD, QRCP, R-SVD, CoR-UTV, RRR-UZVD},
legend style={at={(1,1)},anchor=north east,draw=black,fill=white,legend cell align=left,font=\tiny},
title = {The devil's stairs}
]

\addplot+[smooth,color=red,densely dotted, every mark/.append style={solid}, mark=|]
table[row sep=crcr]{
1	1.00000000000001\\
2	1.00000000000000\\
3	1.00000000000000\\
4	1.00000000000000\\
5	0.999999999999999\\
6	0.999999999999999\\
7	0.999999999999999\\
8	0.999999999999998\\
9	0.999999999999995\\
10	0.999999999999993\\
11	0.251188643150959\\
12	0.251188643150959\\
13	0.251188643150959\\
14	0.251188643150958\\
15	0.251188643150958\\
16	0.251188643150958\\
17	0.251188643150957\\
18	0.251188643150957\\
19	0.251188643150957\\
20	0.251188643150957\\
21	0.0630957344480198\\
22	0.0630957344480195\\
23	0.0630957344480195\\
24	0.0630957344480194\\
25	0.0630957344480194\\
26	0.0630957344480193\\
27	0.0630957344480193\\
28	0.0630957344480193\\
29	0.0630957344480193\\
30	0.0630957344480192\\
31	0.0158489319246112\\
32	0.0158489319246111\\
33	0.0158489319246111\\
34	0.0158489319246111\\
35	0.0158489319246111\\
36	0.0158489319246111\\
37	0.0158489319246111\\
38	0.0158489319246111\\
39	0.0158489319246111\\
40	0.0158489319246111\\
};
\addplot+[smooth,color=black,densely dotted, every mark/.append style={solid}, mark=triangle]
table[row sep=crcr]{
1	0.164065525361425\\
2	0.157950001112986\\
3	0.153095293371917\\
4	0.146517372719229\\
5	0.140238704051759\\
6	0.133011753284716\\
7	0.128060186415024\\
8	0.122900457754931\\
9	0.108366248086257\\
10	0.100369956457925\\
11	0.0771056834619969\\
12	0.0568957539942523\\
13	0.0556046836231045\\
14	0.0521324603391850\\
15	0.0494801242150909\\
16	0.0434296120515032\\
17	0.0405630972984547\\
18	0.0374444742700987\\
19	0.0312737483888488\\
20	0.0289959832432638\\
21	0.0213493043647498\\
22	0.0175462146773521\\
23	0.0146593677840351\\
24	0.0132847828964891\\
25	0.0120208085398750\\
26	0.0109078457770764\\
27	0.0105422058137924\\
28	0.00981581266778682\\
29	0.00880336675701559\\
30	0.00836171742136564\\
31	0.00550823401413481\\
32	0.00498159872419664\\
33	0.00449760106621851\\
34	0.00408763317685733\\
35	0.00400661072884057\\
36	0.00357808709659297\\
37	0.00301170614598169\\
38	0.00271113666893662\\
39	0.00246182621877657\\
40	0.00233831792095783\\
};

\addplot+[smooth,color=olive,densely dotted, every mark/.append style={solid}, mark=diamond]
  table[row sep=crcr]{
  1	1.00000000000000\\
2	0.999999999999998\\
3	0.999999999999997\\
4	0.999999999999997\\
5	0.999999999999996\\
6	0.999999999999996\\
7	0.999999999999994\\
8	0.999999999999993\\
9	0.999999999999977\\
10	0.999999999999967\\
11	0.251188643150922\\
12	0.251188643150788\\
13	0.251188643150514\\
14	0.251188643150084\\
15	0.251188643149366\\
16	0.251188643149095\\
17	0.251188643148344\\
18	0.251188643145334\\
19	0.251188643127815\\
20	0.251188643122955\\
21	0.0630957344118526\\
22	0.0630957342865163\\
23	0.0630957337995500\\
24	0.0630957333755196\\
25	0.0630957332187230\\
26	0.0630957306051668\\
27	0.0630957297585963\\
28	0.0630957201448058\\
29	0.0630957100174593\\
30	0.0630956734831746\\
31	0.0158489209688371\\
32	0.0158486124712117\\
33	0.0158482682694501\\
34	0.0158476312232751\\
35	0.0158474324938597\\
36	0.0158461359438284\\
37	0.0158432755920466\\
38	0.0158384922590068\\
39	0.0158370865948389\\
40	0.0158075382894710\\
};

\addplot+[smooth,color=purple,densely dotted, every mark/.append style={solid}, mark=star]
  table[row sep=crcr]{
  1	1.00000000000000\\
2	0.999999999999998\\
3	0.999999999999997\\
4	0.999999999999997\\
5	0.999999999999996\\
6	0.999999999999996\\
7	0.999999999999994\\
8	0.999999999999993\\
9	0.999999999999977\\
10	0.999999999999967\\
11	0.251188643150922\\
12	0.251188643150788\\
13	0.251188643150514\\
14	0.251188643150084\\
15	0.251188643149366\\
16	0.251188643149095\\
17	0.251188643148344\\
18	0.251188643145334\\
19	0.251188643127815\\
20	0.251188643122955\\
21	0.0630957344118526\\
22	0.0630957342865163\\
23	0.0630957337995500\\
24	0.0630957333755196\\
25	0.0630957332187230\\
26	0.0630957306051668\\
27	0.0630957297585963\\
28	0.0630957201448058\\
29	0.0630957100174593\\
30	0.0630956734831746\\
31	0.0158489209688371\\
32	0.0158486124712117\\
33	0.0158482682694501\\
34	0.0158476312232751\\
35	0.0158474324938597\\
36	0.0158461359438284\\
37	0.0158432755920466\\
38	0.0158384922590068\\
39	0.0158370865948389\\
40	0.0158075382894710\\
};

\addplot+[smooth,color=blue,densely dotted, every mark/.append style={solid}, mark=x]
  table[row sep=crcr]{
  1	0.999934190846255\\
2	0.999925759889860\\
3	0.999911559171821\\
4	0.999757172334876\\
5	0.999740781826369\\
6	0.999505163488621\\
7	0.997271483087945\\
8	0.994547158794357\\
9	0.993025923920002\\
10	0.980992309176834\\
11	0.259300572124159\\
12	0.251463969855461\\
13	0.251457602698757\\
14	0.251274163677450\\
15	0.251164283788927\\
16	0.251156744316190\\
17	0.251125169309715\\
18	0.251085936942232\\
19	0.250700707687012\\
20	0.248930021357893\\
21	0.0636797286253612\\
22	0.0631699238051780\\
23	0.0631407425368120\\
24	0.0631118571330911\\
25	0.0630821665867359\\
26	0.0630607140015002\\
27	0.0626834593522978\\
28	0.0626753598033491\\
29	0.0620616794867621\\
30	0.0532461090434529\\
31	0.0191829739674615\\
32	0.0159625954471212\\
33	0.0158649048579138\\
34	0.0158617903660335\\
35	0.0158507374307986\\
36	0.0158424789996059\\
37	0.0158320681009871\\
38	0.0158227775062661\\
39	0.0158150167461809\\
40	0.0140037494315384\\
  };    
\end{axis}
\end{tikzpicture}%

%% file: GraphsUZV_2018/r_FroErr_Gear_UZV.tex
%
%
%
\usetikzlibrary{positioning,calc}

\definecolor{mycolor1}{rgb}{0.00000,1.00000,1.00000}%
\definecolor{mycolor2}{rgb}{1.00000,0.00000,1.00000}%

\pgfplotsset{every axis label/.append style={font=\footnotesize},
every tick label/.append style={font=\footnotesize}
}

\begin{tikzpicture}[font=\footnotesize] 

\begin{axis}[%
name=ber,
ymode=log,
width  = 0.6\columnwidth,
height = 0.45\columnwidth,
scale only axis,
xmin  = 5,
xmax  = 295,
xlabel= {Approximation rank},
xmajorgrids,
xtick       ={40,120,200,280},
xticklabels ={$40$, $120$ , $200$, $280$},
ymin = 3,
ymax = 199,
ylabel={$\|{\bf A} - \hat{\bf A}_{\text{approx}}\|_F$},
ymajorgrids,
legend entries={Tr SVD, Tr QRCP, R-SVD q=0, R-SVD q=1,UZVD q=0, 
UZVD q=1},
legend style={at={(1,1)},anchor=north east,draw=black,fill=white,legend cell align=left,font=\tiny}
]

\addplot+[smooth,color=gray,densely dotted, every mark/.append style={solid}, mark=asterisk]
table[row sep=crcr]
{
5	113.127160080435\\
15	65.0396321619753\\
25	46.8344672369866\\
35	37.0182723514773\\
45	30.7315037203961\\
55	26.2862146935240\\
65	22.8907671515119\\
75	20.1811640662403\\
85	18.0391136692251\\
95	16.2419808691967\\
105	14.7307471725022\\
115	13.4261407189821\\
125	12.2743649493987\\
135	11.2543825694219\\
145	10.3355346458861\\
155	9.52999472971659\\
165	8.78763799375611\\
175	8.10435605535287\\
185	7.47490437466261\\
195	6.92001490916757\\
205	6.40989416469689\\
215	5.93562749222088\\
225	5.49340294110488\\
235	5.08832432377606\\
245	4.71081467835522\\
255	4.35659021659413\\
265	4.03295568339158\\
275	3.73425114243015\\
285	3.45601877418760\\
295	3.19937756095855 \\
}; 
\addplot+[smooth,color=purple, densely dotted, every mark/.append style={solid}, mark=x]
  table[row sep=crcr]
  {
5	165.589313729307\\
15	103.104852703844\\
25	71.9336368022893\\
35	57.1587303313088\\
45	49.7628174103858\\
55	43.7879245936781\\
65	39.5983070767169\\
75	36.0118591331965\\
85	33.6009925747146\\
95	28.7149053928279\\
105	26.6226271793690\\
115	24.3628013092934\\
125	20.8858967776652\\
135	19.7758408228489\\
145	18.8598745263979\\
155	17.9536066531726\\
165	16.4147514563751\\
175	15.0110650275362\\
185	14.1225816505267\\
195	13.2964592332173\\
205	12.5710527413513\\
215	12.0335421948318\\
225	11.4665250635272\\
235	10.2277131110906\\
245	9.38504414482297\\
255	8.86778023610350\\
265	8.45848602853607\\
275	8.08771285513463\\
285	7.12168948149828\\
295	6.76340442010638\\
}; 
\addplot+[smooth,color=blue, densely dotted, every mark/.append style={solid}, mark=triangle]
  table[row sep=crcr]
  {
5	198.736254322996\\
15	104.174184171521\\
25	82.2920334376797\\
35	61.8684558242737\\
45	52.0593724169386\\
55	45.4509023194245\\
65	39.0721960111490\\
75	35.3058581724190\\
85	31.3859724502080\\
95	28.4634929911408\\
105	26.6255792586465\\
115	24.0668754008060\\
125	22.2411753181908\\
135	20.4723103459044\\
145	19.1761684061879\\
155	17.7415242787572\\
165	16.6567164962503\\
175	15.2388504545478\\
185	14.2029669718233\\
195	13.4573932306739\\
205	12.4914130759158\\
215	11.8249199417053\\
225	10.9774601445606\\
235	10.4759707806572\\
245	9.64618531084318\\
255	9.14847990134760\\
265	8.56047342000471\\
275	8.01278400011652\\
285	7.48283184270495\\
295	6.92528703243961\\
};

\addplot+[smooth,color=olive,densely dotted, every mark/.append style={solid}, mark = diamond]
  table[row sep=crcr]
{
5	119.215268404256\\
15	70.5718495109994\\
25	49.6320822839156\\
35	39.0606549925464\\
45	32.9907389187450\\
55	28.0133227680193\\
65	24.4902534207153\\
75	21.7622857644341\\
85	19.2314289723909\\
95	17.4246887229794\\
105	15.8654520951009\\
115	14.4666911631032\\
125	13.2650428924215\\
135	12.2096512655925\\
145	11.1830353554933\\
155	10.3173476481698\\
165	9.54318085647258\\
175	8.77953904675649\\
185	8.17011417655641\\
195	7.56759192397715\\
205	6.99095011043945\\
215	6.46659070375483\\
225	6.00274760138303\\
235	5.57110551613331\\
245	5.15363196973522\\
255	4.77439520630524\\
265	4.43701142114874\\
275	4.10392764117562\\
285	3.81859133183317\\
295	3.50566920264563 \\
};

\addplot+[smooth,color = black ,densely dotted, every mark/.append style={solid}, mark=star]
  table[row sep=crcr]
  {
5	175.523385804335\\
15	99.6083241091235\\
25	75.2693855279645\\
35	59.3206365604519\\
45	51.6594697068012\\
55	44.8627917921231\\
65	40.4441902393841\\
75	35.0281480622496\\
85	31.6834352247104\\
95	29.0668418483956\\
105	26.3465932032612\\
115	24.4914459437500\\
125	21.9626053289945\\
135	20.7537397395190\\
145	19.0094982579005\\
155	17.5277421444826\\
165	16.5425980149950\\
175	15.3481269967419\\
185	14.4011703913304\\
195	13.3535339899693\\
205	12.4910340974367\\
215	11.6545438173921\\
225	10.9497151771003\\
235	10.2773609186726\\
245	9.65409284807292\\
255	8.95802009640951\\
265	8.37280805251357\\
275	7.92697465597522\\
285	7.46463163923867\\
295	6.90169242292872 \\
};
\addplot+[smooth,color = violet,densely dotted,every mark/.append style={solid},mark=+]
  table[row sep=crcr]
  {
5	120.395005376443\\
15	70.0446549109648\\
25	49.1693824673539\\
35	39.3574762968164\\
45	32.8674233033153\\
55	28.0943614753176\\
65	24.6931054299396\\
75	21.6935947176389\\
85	19.3276241731625\\
95	17.4059268564784\\
105	15.7402692375969\\
115	14.4686766933353\\
125	13.2772049389625\\
135	12.1836387565988\\
145	11.1719896080367\\
155	10.3388841259223\\
165	9.52941600354573\\
175	8.81653754139914\\
185	8.17485565294307\\
195	7.52633694299458\\
205	6.99472270538392\\
215	6.48830119640714\\
225	6.02372712396301\\
235	5.56213357237267\\
245	5.15330494104008\\
255	4.76878640286080\\
265	4.44335427132534\\
275	4.09927098999832\\
285	3.78513379101605\\
295	3.52751382403797 \\
};

\end{axis}

\end{tikzpicture}%

%% file: GraphsUZV_2018/r_Time_Gear_UZV.tex
%
%
%
\usetikzlibrary{positioning,calc}

\definecolor{mycolor1}{rgb}{0.00000,1.00000,1.00000}%
\definecolor{mycolor2}{rgb}{1.00000,0.00000,1.00000}%

\pgfplotsset{every axis label/.append style={font=\footnotesize},
every tick label/.append style={font=\footnotesize}
}

\begin{tikzpicture}[font=\footnotesize] 

\begin{axis}[%
name=ber,
ymode=log,
width  = 0.66\columnwidth,
height = 0.51\columnwidth,
scale only axis,
xmin  = 5,
xmax  = 295,
xlabel= {Approximation rank},
xmajorgrids,
xtick       ={40,120,200,280},
xticklabels ={$40$, $120$ , $200$, $280$},
ymin = 0.0,
ymax = 3.99,
ylabel={Time (seconds)},
ymajorgrids,
]

\addplot+[smooth,color=gray,densely dotted, every mark/.append style={solid}, mark=asterisk]
table[row sep=crcr]
{
5	0.0941121973015007\\
15	0.0879612365282019\\
25	0.108391567720097\\
35	0.178942016956795\\
45	0.218702766170434\\
55	0.272007329561355\\
65	0.333234346850006\\
75	0.408021600036675\\
85	0.488040539754591\\
95	0.582853397140842\\
105	0.689877172371276\\
115	0.802993063875704\\
125	0.932285034543774\\
135	1.03140072488212\\
145	1.15664508350787\\
155	1.29719414374069\\
165	1.43127133540749\\
175	1.56951109116176\\
185	1.76493230487746\\
195	1.89827633513721\\
205	2.10077975803963\\
215	2.29348180849674\\
225	2.41258923323565\\
235	2.60722186223682\\
245	2.80683814136228\\
255	3.03394609159625\\
265	3.21485187265777\\
275	3.54704241388347\\
285	3.70054239575115\\
295	3.98009961142103 \\
}; 

\addplot+[smooth,color=blue, densely dotted, every mark/.append style={solid}, mark=triangle]
  table[row sep=crcr]
  {
5	0.00852937332660955\\
15	0.00928374612293437\\
25	0.0119074581297439\\
35	0.0129478615156290\\
45	0.0153294399714262\\
55	0.0159208582538647\\
65	0.0155890026473182\\
75	0.0176831141603814\\
85	0.0197878313680867\\
95	0.0225459962134249\\
105	0.0240325040592441\\
115	0.0269428140944892\\
125	0.0289518357773187\\
135	0.0308721507459225\\
145	0.0325362184472031\\
155	0.0358106411381895\\
165	0.0390552994603412\\
175	0.0418723088049845\\
185	0.0435021645880646\\
195	0.0470779942783988\\
205	0.0513537998151188\\
215	0.0541933890257746\\
225	0.0570275043295182\\
235	0.0612954411250521\\
245	0.0675900919547937\\
255	0.0721371980028451\\
265	0.0716986012115125\\
275	0.0758016365613191\\
285	0.0775820247844820\\
295	0.0833604177685755\\
}; 
\addplot+[smooth,color=olive,densely dotted, every mark/.append style={solid}, mark = diamond]
  table[row sep=crcr]
  {
5	0.0101373334820411\\
15	0.0112334833411907\\
25	0.0130627946073806\\
35	0.0152192436565972\\
45	0.0162674249852376\\
55	0.0184344078840681\\
65	0.0191911755146670\\
75	0.0226674485230373\\
85	0.0240951118695513\\
95	0.0275176721663466\\
105	0.0292392158902044\\
115	0.0322912611128863\\
125	0.0354493632819904\\
135	0.0374658137407380\\
145	0.0397751182192833\\
155	0.0433068145351387\\
165	0.0467382699306661\\
175	0.0485025785522748\\
185	0.0533216693500215\\
195	0.0550685298933478\\
205	0.0586898614348889\\
215	0.0624384613121365\\
225	0.0657994401561706\\
235	0.0736452593571307\\
245	0.0771995355389985\\
255	0.0807777600636068\\
265	0.0841955306918539\\
275	0.0862273765137919\\
285	0.0914997752276974\\
295	0.0950927108771320 \\
};

\addplot+[smooth,color = black ,densely dotted, every mark/.append style={solid}, mark=star]
  table[row sep=crcr]
{
5	0.0106525649701433\\
15	0.0110323172621707\\
25	0.0124271371672121\\
35	0.0142037620793730\\
45	0.0149495818961477\\
55	0.0166252816496166\\
65	0.0177063782647578\\
75	0.0203762763611383\\
85	0.0228730621514233\\
95	0.0262035923882587\\
105	0.0275248566691687\\
115	0.0301872568815563\\
125	0.0319231408730471\\
135	0.0348944459687755\\
145	0.0381216562126449\\
155	0.0405127271756898\\
165	0.0459753872618081\\
175	0.0471926042043711\\
185	0.0516760760845691\\
195	0.0537623188564461\\
205	0.0580966267732893\\
215	0.0618031459911501\\
225	0.0634182906494038\\
235	0.0691600768810226\\
245	0.0717475242545394\\
255	0.0759381421149397\\
265	0.0852684164466270\\
275	0.0858455715066724\\
285	0.094407276464458\\
295	0.0993446720753429\\
};
\addplot+[smooth,color = violet ,densely dotted, every mark/.append style={solid}, mark=+]
  table[row sep=crcr]
{
5	0.0111103204356683\\
15	0.0118653774703571\\
25	0.0133258842583438\\
35	0.0157634834301417\\
45	0.0171059591003360\\
55	0.0188384506380181\\
65	0.0215644562802476\\
75	0.0233178170880321\\
85	0.0251933144437929\\
95	0.0280113501459823\\
105	0.0322594440289597\\
115	0.0331951399917481\\
125	0.0366132527391772\\
135	0.0389926916500340\\
145	0.0410313798556120\\
155	0.0483954952483067\\
165	0.0490017304388226\\
175	0.0524824509965590\\
185	0.0569525802286588\\
195	0.0588239721542355\\
205	0.0633987898560294\\
215	0.0690307558302240\\
225	0.0709346490780914\\
235	0.0773900959233763\\
245	0.0796364504724324\\
255	0.0874607161649261\\
265	0.0899284218247406\\
275	0.0955347288602847\\
285	0.096152398563373\\
295	0.0965501281236337\\
};
\end{axis}

\end{tikzpicture}%